# Mathematical base of Difference Operator of Particle Method

H. Isshiki, Institute of Mathematical Analysis, Osaka, Japan, isshiki@dab.hi-ho.ne.jp

D. Kitazawa, Institute of Industrial Science, The University of Tokyo, dkita@iis.u-tokyo.ac.jp

**Abstract**
Mathematical base of difference operators in Moving Particle Semi-implicit method (MPS) are not given sufficiently in contrast to Smooth Particle Hydrodynamics method (SPH). Iribe and Nakaza proposed a method to improve the accuracy of the gradient operator, and Khayyer and Gotoh gave an ingenuity also for gradient operator too. An extension to higher order difference operators of Iribe-Nakaza method is given in this paper. The proposed method is a special case of the author's method called Discrete Differential Operators on Irregular Nodes (DDIN).

## 1. Introduction

Moving Particle Semi-implicit method (MPS)[1] is a particle method and is widely used. It gives plausible numerical results in many cases. However, mathematical base of difference operators in Moving Particle Semi-implicit method (MPS) are not given sufficiently in contrast to Smooth Particle Hydrodynamics method (SPH) based on approximations of Dirac's delta function. Iribe and Nakaza[2,3] proposed a method to improve the accuracy of the gradient operator, and Khayyer and Gotoh[4,5] gave an ingenuity also for gradient operator too.

An extension to higher order difference operators of the former method is given in this paper. The proposed method in the present paper uses Taylor expansion. The method is a special case of the author's method called Discrete Differential operators on Irregular Nodes (DDIN)[6]. In the author's paper[6], interpolation of discrete data given on irregular mesh is used. If we use the interpolation using power functions, the author's method can give a similar result obtained by Iribe and Nakaza for the gradient operator.

## 2. Summary of MPS method

Let values of a function $\phi$ be given on irregular nodes $(x_j, y_j)$, $j=1,2,\cdots,J$ as shown in Fig. 1. In Koshizuka's original MPS method[1], the gradient and Laplace operators at a point $C(x_i, y_i)$ are given as

$$\langle \nabla \phi \rangle_i = \frac{d}{n^0} \sum_{j \neq i} \frac{\phi_j - \phi_i}{|\mathbf{r}_j - \mathbf{r}_i|^2} (\mathbf{r}_j - \mathbf{r}_i) w(|\mathbf{r}_j - \mathbf{r}_i|), \qquad (1)$$

$$\langle \nabla^2 \phi \rangle_i = \frac{2d}{\lambda n^0} \sum_{j \neq i} (\phi_j - \phi_i) w(|\mathbf{r}_j - \mathbf{r}_i|), \qquad (2)$$

where $d$ and $n^0$ are the number of dimension and constant particle number density, respectively, and $w$ is a weight function defined as

$$w(r) = \begin{cases} r_e/r - 1 & r < r_e \\ 0 & r_e < r \end{cases}. \qquad (3)$$

$\langle n \rangle_i$ is the particle number density at $C(x_i, y_i)$ defined as

$$\langle n \rangle_i = \sum_{j \neq i} w(|\mathbf{r}_j - \mathbf{r}_i|), \qquad (4)$$

and $\lambda$ is a constant defined as

$$\lambda = \frac{\sum_{j \neq i} |\mathbf{r}_j - \mathbf{r}_i|^2 w(|\mathbf{r}_j - \mathbf{r}_i|)}{\sum_{j \neq i} w(|\mathbf{r}_j - \mathbf{r}_i|)}. \tag{5}$$

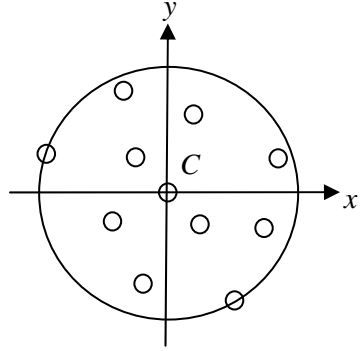

Figure 1. Definition of neighboring nodes surrounding node $C$

The mathematical background of Eq. (1) might be given as follows. If $\langle \nabla \phi \rangle_i$ is parallel to $(\mathbf{r}_j - \mathbf{r}_i)$, then we have

$$\langle \nabla \phi \rangle_i = \frac{\phi_j - \phi_i}{|\mathbf{r}_j - \mathbf{r}_i|} \frac{\mathbf{r}_j - \mathbf{r}_i}{|\mathbf{r}_j - \mathbf{r}_i|}. \tag{6}$$

Equation (6) could suggest Eq. (1).

The mathematical background of Eq. (2) might be given as follows. In a regular equi-spacing mesh, we have

$$\frac{d^2 \phi}{dx^2} = \frac{\phi_{i+1} - 2\phi_i + \phi_{i-1}}{\Delta^2} = \frac{(\phi_{i+1} - \phi_i) + (\phi_{i-1} - \phi_i)}{\Delta^2} \quad \text{in 1d}, \tag{7}$$

and

$$\begin{aligned} \frac{\partial^2 \phi}{\partial x^2} + \frac{\partial^2 \phi}{\partial y^2} &= \frac{\phi_{i+1,j} - 2\phi_{i,j} + \phi_{i-1,j}}{\Delta^2} + \frac{\phi_{i,j+1} - 2\phi_{i,j} + \phi_{i,j-1}}{\Delta^2} \\ &= \frac{(\phi_{i+1,j} - \phi_{i,j}) + (\phi_{i-1,j} - \phi_{i,j}) + (\phi_{i,j+1} - \phi_{i,j}) + (\phi_{i,j-1} - \phi_{i,j})}{\Delta^2} \end{aligned} \quad \text{in 2d}. \tag{8}$$

Equations (7) and (8) could suggest Eq. (2).

## 3. Summary of the method to improve the accuracy of difference operators
### 3.1. Iribe-Nakaza's method to improve gradient operator

Iribe and Nakaza proposed a method[2,3] to improve the accuracy of the gradient operator given by Eq. (1). The first order approximation of a function $\phi$ is given by

$$\phi_j = \phi_i + (\nabla \phi)_i \cdot (\mathbf{r}_j - \mathbf{r}_i), \tag{9}$$

where $j = 1, 2, \cdots, J$. Rewriting, we have

$$(\nabla \phi)_i \cdot (\mathbf{r}_j - \mathbf{r}_i) = \phi_j - \phi_i. \tag{10}$$

Dividing the both sides by $|\mathbf{r}_j - \mathbf{r}_i|$, we obtain

$$(\nabla \phi)_i \cdot \frac{(\mathbf{r}_j - \mathbf{r}_i)}{|\mathbf{r}_j - \mathbf{r}_i|} = \frac{\phi_j - \phi_i}{|\mathbf{r}_j - \mathbf{r}_i|}. \tag{11}$$

Rewriting again with matrix notation, we have

$$\begin{bmatrix} \dfrac{(\mathbf{r}_j-\mathbf{r}_i)_1}{|\mathbf{r}_j-\mathbf{r}_i|} & \dfrac{(\mathbf{r}_j-\mathbf{r}_i)_2}{|\mathbf{r}_j-\mathbf{r}_i|} & \cdots & \dfrac{(\mathbf{r}_j-\mathbf{r}_i)_N}{|\mathbf{r}_j-\mathbf{r}_i|} \end{bmatrix} \begin{bmatrix} (\nabla\phi)_{i1} \\ (\nabla\phi)_{i2} \\ \vdots \\ (\nabla\phi)_{iN} \end{bmatrix} = \dfrac{\phi_j-\phi_i}{|\mathbf{r}_j-\mathbf{r}_i|}, \quad (12)$$

where $N$ is the number of dimension. Changing $j=1,2,\cdots,J$, we finally obtain

$$\begin{bmatrix} \dfrac{(\mathbf{r}_1-\mathbf{r}_i)_1}{|\mathbf{r}_1-\mathbf{r}_i|} & \dfrac{(\mathbf{r}_1-\mathbf{r}_i)_2}{|\mathbf{r}_1-\mathbf{r}_i|} & \cdots & \dfrac{(\mathbf{r}_1-\mathbf{r}_i)_N}{|\mathbf{r}_1-\mathbf{r}_i|} \\ \dfrac{(\mathbf{r}_2-\mathbf{r}_i)_1}{|\mathbf{r}_2-\mathbf{r}_i|} & \dfrac{(\mathbf{r}_2-\mathbf{r}_i)_2}{|\mathbf{r}_2-\mathbf{r}_i|} & \cdots & \dfrac{(\mathbf{r}_2-\mathbf{r}_i)_2}{|\mathbf{r}_2-\mathbf{r}_i|} \\ \vdots & \vdots & \vdots & \vdots \\ \dfrac{(\mathbf{r}_J-\mathbf{r}_i)_1}{|\mathbf{r}_J-\mathbf{r}_i|} & \dfrac{(\mathbf{r}_J-\mathbf{r}_i)_2}{|\mathbf{r}_J-\mathbf{r}_i|} & \cdots & \dfrac{(\mathbf{r}_J-\mathbf{r}_i)_N}{|\mathbf{r}_J-\mathbf{r}_i|} \end{bmatrix} \begin{bmatrix} (\nabla\phi)_{i1} \\ (\nabla\phi)_{i2} \\ \vdots \\ (\nabla\phi)_{iN} \end{bmatrix} = \begin{bmatrix} \dfrac{\phi_1-\phi_i}{|\mathbf{r}_1-\mathbf{r}_i|} \\ \dfrac{\phi_2-\phi_i}{|\mathbf{r}_2-\mathbf{r}_i|} \\ \vdots \\ \dfrac{\phi_J-\phi_i}{|\mathbf{r}_J-\mathbf{r}_i|} \end{bmatrix} \quad (13)$$

or

$$[A][\nabla\phi] = [b]. \quad (14)$$

where

$$[A]_{jk} = \dfrac{(\mathbf{r}_j-\mathbf{r}_i)_k}{|\mathbf{r}_j-\mathbf{r}_i|}, \quad [\nabla\phi]_k = (\nabla\phi)_k, \quad [b]_j = \dfrac{\phi_j-\phi_i}{|\mathbf{r}_j-\mathbf{r}_i|}. \quad (15)$$

Equation (15) can be solved by Least Square method (LSM). Multiplying $[A]^T[W]$ on both sides, we derive

$$[A]^T[W][A][\nabla\phi] = [A]^T[W][b], \quad (16)$$

where $[W]$ is weight function. Finally, we obtain

$$[\nabla\phi] = \left([A]^T[W][A]\right)^{-1}[A]^T[W][b]. \quad (17)$$

Although $|r_j - r_i|$ is used instead of $|\mathbf{r}_j - \mathbf{r}_i|$ in Ref. (2), the difference is not critical, namely, the result does not change much. Equations (16) and (17) corresponds to Eqs. (16), (17) and (18) in Ref. (2). Equation (16) is identical to Eq. (15) in Ref. (2). Equation (17) could be more reasonable than those in (16) and (17) in Ref. (2).

**3.2. Khayyer and Gotoh's method to improve gradient operator**

Khayyer and Gotoh proposed a method[4] also to improve the accuracy of the gradient operator for pressure. MPS is invented to solve Navier-Stokes equation for viscous fluid. Gradient of pressure $p$ appears in the equation:

$$\langle p \rangle_i = \dfrac{D_S}{n_0} \sum_{j \neq i} \dfrac{p_j - \hat{p}_i}{|\mathbf{r}_j-\mathbf{r}_i|^2}(\mathbf{r}_j-\mathbf{r}_i)w(|\mathbf{r}_j-\mathbf{r}_i|), \quad (18)$$

$$\hat{p}_i = \min_{j \in J}(p_i, p_j), \quad J = \{j : w(|\mathbf{r}_j-\mathbf{r}_i|) \neq 0\}, \quad (19)$$

where $D_S$ and $n_0$ are the number of dimension and constant particle density, respectively. $\hat{p}_i$ means the minimum pressure of the related particles.

Khayyer and Gotoh noticed that the pressure gradient force from particle $j$ to $i$:

$$\mathbf{A}^p_{j \to i} = -\dfrac{mD_S}{\rho n_0} \sum_{j \neq i} \dfrac{p_j - \hat{p}_i}{|\mathbf{r}_j-\mathbf{r}_i|^2}(\mathbf{r}_j-\mathbf{r}_i)w(|\mathbf{r}_j-\mathbf{r}_i|) \quad (20)$$

is not equal to the pressure gradient force from particle $i$ to $j$:

$$\mathbf{A}^p_{i \to j} = -\dfrac{mD_S}{\rho n_0} \sum_{j \neq i} \dfrac{p_i - \hat{p}_j}{|\mathbf{r}_i-\mathbf{r}_j|^2}(\mathbf{r}_i-\mathbf{r}_j)w(|\mathbf{r}_i-\mathbf{r}_j|). \quad (21)$$

Namely

$$\mathbf{A}^p_{j \to i} \neq -\mathbf{A}^p_{i \to j} . \quad (22)$$

Equation (22) is against Newton's third law and breaks momentum conservation law. Khayyer-Gotoh proposes an ingenuity to replace Eq. (18) with

$$\langle p \rangle_i = \frac{D_S}{n_0} \sum_{j \neq i} \frac{(p_i + p_j) - (\hat{p}_i + \hat{p}_j)}{|\mathbf{r}_j - \mathbf{r}_i|^2} (\mathbf{r}_j - \mathbf{r}_i) w(|\mathbf{r}_j - \mathbf{r}_i|) . \quad (23)$$

## 4. New proposals to improve difference operators in irregular mesh
### 4.1. Method using Taylor expansion

In order to obtain the higher order difference formula, we must extend Iribe-Nakaza's theory. For simplicity, we consider 2d case.

We notice that Eq. (10) is the first three terms of Taylor expansion. Rewriting Eq. (9), we have

$$\phi_j = \phi_i + (\phi_x)_i (x_j - x_i) + (\phi_y)_i (y_j - y_i) . \quad (24)$$

Rewriting Eq. (17), we obtain

$$\begin{bmatrix} x_1 - x_i & y_1 - y_i \\ x_2 - x_i & y_2 - y_i \\ & \vdots \\ x_J - x_i & y_J - y_i \end{bmatrix} \begin{bmatrix} (\phi_x)_i \\ (\phi_y)_i \end{bmatrix} = \begin{bmatrix} \phi_1 - \phi_i \\ \phi_2 - \phi_i \\ \vdots \\ \phi_J - \phi_i \end{bmatrix} . \quad (25)$$

The solution would be obtained by LSM.

If we use the higher Taylor expansion, we could derive a formula for the higher order difference operator. If we take the terms whose order are lower than or equal to the second, we have

$$\phi_j = \phi_i + (\phi_x)_i (x_j - x_i) + (\phi_y)_i (y_j - y_i) + \frac{1}{2}(\phi_{xx})_i (x_j - x_i)^2 + (\phi_{xx})_i (x_j - x_i)(y_j - y_i) + \frac{1}{2}(\phi_{yy})_i (y_j - y_i)^2 . \quad (26)$$

Rewriting Eq. (26), we obtain

$$\begin{bmatrix} x_1 - x_i & y_1 - y_i & (x_1 - x_i)^2 & (x_1 - x_i)(y_1 - y_i) & (y_1 - y_i)^2 \\ x_2 - x_i & y_2 - y_i & (x_2 - x_i)^2 & (x_2 - x_i)(y_2 - y_i) & (y_2 - y_i)^2 \\ \vdots & \vdots & \vdots & \vdots & \vdots \\ \vdots & \vdots & \vdots & \vdots & \vdots \\ \vdots & \vdots & \vdots & \vdots & \vdots \\ x_J - x_i & y_J - y_i & (x_J - x_i)^2 & (x_J - x_i)(y_J - y_i) & (y_N - y_i)^2 \end{bmatrix} \begin{bmatrix} (\phi_x)_i \\ (\phi_y)_i \\ (\phi_{xx})_i \\ (\phi_{xy})_i \\ (\phi_{yy})_i \end{bmatrix} = \begin{bmatrix} \phi_1 - \phi_i \\ \phi_2 - \phi_i \\ \vdots \\ \vdots \\ \phi_J - \phi_i \end{bmatrix} . \quad (27)$$

The solution would be obtained by LSM.

### 4.3. Method using interpolation of discrete data at irregular nodes

Taylor series given by Eqs. (24) and (26) may be considered as an interpolation of a function $\phi(\mathbf{x})$ using power functions. If we use functions $\Phi_\mu(\mathbf{x})$, $\mu = 1, 2, \cdots, M$, $\phi(\mathbf{x})$ could be interpolated as

$$\phi(\mathbf{x}) = \phi_i + a_1 \Phi_1(\mathbf{x} - \mathbf{x}_i) + a_2 \Phi_2(\mathbf{x} - \mathbf{x}_i) + \cdots + a_M \Phi_M(\mathbf{x} - \mathbf{x}_i) , \quad (28)$$

where $\phi_i = \phi(\mathbf{x}_i)$. The coefficients $a_i$, $i = 1, 2, \cdots, M$ are determined by solving

$$\phi_j - \phi_i = a_1 \Phi_1(\mathbf{x}_j - \mathbf{x}_i) + a_2 \Phi_2(\mathbf{x}_j - \mathbf{x}_i) + \cdots + a_M \Phi_M(\mathbf{x}_j - \mathbf{x}_i) , \quad j = 1, 2, \cdots, J \quad (29)$$

or

$$\begin{bmatrix} \Phi_1(\mathbf{x}_1 - \mathbf{x}_i) & \Phi_2(\mathbf{x}_1 - \mathbf{x}_i) & \cdots & \cdots & \Phi_M(\mathbf{x}_1 - \mathbf{x}_i) \\ \Phi_1(\mathbf{x}_2 - \mathbf{x}_i) & \Phi_2(\mathbf{x}_2 - \mathbf{x}_i) & \cdots & \cdots & \Phi_M(\mathbf{x}_2 - \mathbf{x}_i) \\ \vdots & \vdots & \vdots & \vdots & \vdots \\ \vdots & \vdots & \vdots & \vdots & \vdots \\ \vdots & \vdots & \vdots & \vdots & \vdots \\ \Phi_1(\mathbf{x}_J - \mathbf{x}_i) & \Phi_2(\mathbf{x}_J - \mathbf{x}_i) & \cdots & \cdots & \Phi_M(\mathbf{x}_J - \mathbf{x}_i) \end{bmatrix} \begin{bmatrix} a_1 \\ a_1 \\ \vdots \\ \vdots \\ \vdots \\ a_M \end{bmatrix} = \begin{bmatrix} \phi_1 - \phi_i \\ \phi_2 - \phi_i \\ \vdots \\ \vdots \\ \vdots \\ \phi_J - \phi_i \end{bmatrix} \quad (30)$$

using LSM.

Let $L$ be a linear differential or integral operator. $(L\phi)_i$ is given by

$$(L\phi)_i = a_1 (L\Phi_1(\mathbf{x} - \mathbf{x}_i))_i + a_2 (L\Phi_2(\mathbf{x} - \mathbf{x}_i))_i + \cdots + a_M (L\Phi_M(\mathbf{x} - \mathbf{x}_i))_i. \quad (31)$$

The author has already published a paper on this method[6]. The method is called Discrete Differential operators on Irregular Nodes method (DDIN)[6]. If we use the interpolation using power functions, DDIN can give a similar result obtained by Iribe and Nakaza for the gradient operator.

We apply the method for obtaining the numerical solution of the following boundary value problem :

$$\frac{d^2\phi}{dx^2} = x, \quad -1 < x < 1, \quad (32)$$

$$\phi(-1) = \phi(1) = 0. \quad (33)$$

The exact solution $\phi_{exact}$ is given by

$$\phi_{exact} = \frac{1}{6}x(x^2 - 1). \quad (34)$$

The results are shown in Fig. 2. $M$ is the number of the local neighboring nodes. If the method is applied to regular nodes when $M = 3$, it becomes the ordinary finite difference method. Fig. 2(a) shows the results. If the method is applied to irregular nodes, we obtain the results as shown in Figs. 2(b) and 2(c). The regular nodes $x_i$, $i = 0,1,\cdots,M_{total}$ are generated by $x_i = -1 + idx$, where $dx = 2/(M_{total} - 1)$. The random nodes $x_{Ri}$ are made by adding the random number of uniform distribution $|\delta x_i| < \rho_{rnd} dx$ to $x_i$: $x_{Ri} = x_i + \delta x_i$. The numerical results were correct both for the regular and irregular nodes.

The figures show the results for seven different distributions of particles. The exact solution and the mean of the seven distributions are also shown in the figures.

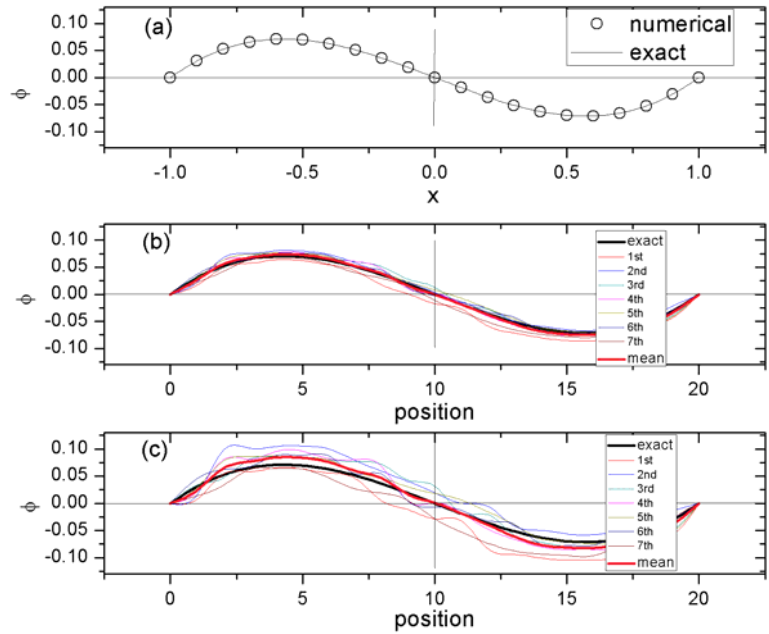

Fig 2. Solution of one dimensional differential equation by MPS ($M = 3$, $dx = 0.1$); (a) Regular nodes

($\rho_{rnd} = 0$); (b) Irregular nodes ($\rho_{rnd} = 0.25$); (c) Irregular nodes ($\rho_{rnd} = 0.5$).

## 5. Conclusions

Mathematical base of difference operators in Moving Particle Semi-implicit method (MPS) are not given sufficiently in contrast to Smooth Particle Hydrodynamics method (SPH) based on approximations of Dirac's delta function. Iribe and Nakaza[2,3] proposed a method to improve the accuracy of the gradient operator, and Khayyer and Gotoh[4] gave an ingenuity also for Laplace operator too.

An extension to higher order difference operators of the Iribe-Nakaza's method is discussed in this paper. The proposed method in the present paper uses Taylor expansion. The method is a special case of the author's method called Discrete Differential operators on Irregular Nodes (DDIN)[6]. In the author's paper[5], interpolation of discrete data given on irregular mesh is used. If we use the interpolation using power functions, the author's method can give a similar result obtained by Iribe and Nakaza for the gradient operator.